\documentclass[11pt]{article}
\usepackage{latexsym,amssymb}
\usepackage[centertags]{amsmath}
\usepackage{fullpage}
\usepackage{fancyhdr}
\title{Convergence of the self-dual Ginzburg-Landau 
gradient flow
}
\author{Sophia Demoulini
\\{\it{\small Centre for Mathematical Sciences, Wilberforce Road, 
Cambridge, CB3 OWB,
England}}\\{\it\small{ email:sd290@cam.ac.uk}}}
\date{} 
\begin{document}

\newtheorem{lemma}{Lemma}[section]
\newtheorem{theorem}[lemma]{Theorem}
\newtheorem{maintheorem}[lemma]{Main Theorem}
\newtheorem{corl}[lemma]{Corollary}
\newtheorem{definition}[lemma]{Definition}
\newtheorem{remark}[lemma]{Remark}
\newtheorem{Rem}{Remark}
\newtheorem{Notation}[lemma]{Notation}
\newcommand{\proof}{\noindent {\it Proof}\;\;\;}
\newcommand{\qed}{\protect~\protect\hfill $\Box$}

\newcommand{\be}{\begin{equation}}
\newcommand{\ee}{\end{equation}}
\newcommand{\ba}{\begin{eqnarray}}
\newcommand{\ea}{\end{eqnarray}}
\newcommand{\bes}{\[}
\newcommand{\ees}{\]}
\newcommand{\bas}{\begin{eqnarray*}}
\newcommand{\eas}{\end{eqnarray*}}
\newcommand{\hoa}{{H^{1}_{\A}}}
\newcommand{\hta}{{H^{2}_{\mbA}}}
\newcommand{\linf}{{L^\infty}}
\newcommand{\tp}{{\tilde P}}
\newcommand{\cf}{{\cal F}}
\newcommand{\ch}{{\cal H}}
\newcommand{\cfnd}{{\cal F}^{n}-{\cal F}^{n-1}}

\newcommand{\spsi}{_{{{\Psi}}}}
\newcommand{\pt}{\frac{\partial}{\partial t}}
\newcommand{\pxk}{\frac{\partial}{\partial {x^k}}}
\newcommand{\pxj}{\frac{\partial}{\partial {x^j}}}
\renewcommand{\theequation}{\arabic{equation}}
\newcommand{\rgt}{\rightarrow}
\newcommand{\lngrgt}{\longrightarrow}
\newcommand{\intsT}{ \int_{0}^{T}\!\!\int_{\Sigma} }
\newcommand{\dxdt}{\;dx\,dt}
\newcommand{\sublt}{_{L^2}}
\newcommand{\sublf}{_{L^4}}
\newcommand{\naf}{\nabla_\mbA \Phi}
\newcommand{\covt}{(\pt -iA_0) }
\newcommand{\ano}{A^{n}_{0}}
\newcommand{\aNo}{A^{N}_{0}}
\newcommand{\Psino}{\Psi^{n}_{0}}
\newcommand{\PsiNo}{\Psi^{N}_{0}}
\newcommand{\Nmo}{{N\!-\!1}}
\newcommand{\nmo}{{n-1}}
\newcommand{\nmt}{{n-2}}
\newcommand{\Nmt}{{N\!-\!2}}
\newcommand{\gotm}{\frac{\gamma}{2\mu}}
\newcommand{\ootm}{\frac{1}{2\mu}}
\newcommand{\tloc}{T_{loc}}
\newcommand{\tmax}{T_{max}}
\font\msym=msbm10
\def\Real{{\mathop{\hbox{\msym \char '122}}}}
\font\smallmsym=msbm7
\def\smr{{\mathop{\hbox{\smallmsym \char '122}}}}
\def\Complex{{\mathop{\hbox{\msym\char'103}}}}
\newcommand{\wkarr}{\; \rightharpoonup \;}
\def\Weak{\,\,\relbar\joinrel\rightharpoonup\,\,}
\newcommand{\To}{\longrightarrow}
\newcommand{\pa}{\partial_{\A}}
\newcommand{\pbfa}{\partial_{\mathbf A}}
\newcommand{\pao}{\partial_{A_1}}
\newcommand{\pat}{\partial_{A_2}}
\newcommand{\dbar}{\bar{\partial}}
\newcommand{\barpa}{\bar{\partial}_{\A}}
\newcommand{\barpbfa}{\bar{\partial}_{\mathbf A}}
\newcommand{\barpaphi}{\bar{\partial}_{\A}\Phi}
\newcommand{\barpbfaphi}{\bar{\partial}_{\mathbf A}\Phi}
\newcommand{\myqed}{\hfill $\Box$}
\newcommand{\cd}{{\cal D}}
\newcommand{\dt}{\hbox{det}\,}
\newcommand{\sma}{_{{ A}}}
\newcommand{\bfpi}{{\mbox{\boldmath$\pi$}}}
\newcommand{\ce}{{\cal E}}
\newcommand{\ulh}{{\underline h}}
\newcommand{\ulg}{{\underline g}}
\newcommand{\ulX}{{\underline {\bf X}}}
\newcommand\la{\label}
\newcommand{\lamo}{\stackrel{\circ}{\lambda}}
\newcommand{\bfjo}{\underline{{\bf J}}}
\newcommand{\Vflato}{V^\flat_0}
\newcommand{\cm}{{\cal M}}
\newcommand{\dist}{{\mbox{dist}}}
\newcommand{\cs}{{\cal S}}
\newcommand{\mcF}{{\mathcal F}}
\newcommand{\mcn}{{\mathcal  V}}
\newcommand{\mce}{{\mathcal  E}}
\newcommand{\mcb}{{\mathcal B}}
\newcommand{\mca}{{\mathcal A}}
\newcommand{\mcdpsi}{{\mathcal D}_{\psi}}
\newcommand{\mcd}{{\mathcal D}}
\newcommand{\mcl}{{\mathcal L}}
\newcommand{\mclaphi}{{\mathcal L}_{(\mbA,\Phi)}}
\newcommand{\mcdadjp}{\mcdadj_{\psi}}
\newcommand{\mcdadj}{{\mathcal D}^{\ast}}
\newcommand{\zl}{{Z_\Lambda}}
\newcommand{\thl}{{\Theta_\Lambda}}
\newcommand{\ca}{{\cal A}}
\newcommand{\cb}{{\cal B}}
\newcommand{\cg}{{\cal G}}
\newcommand{\cu}{{\cal U}}
\newcommand{\cl}{{\cal L}}
\newcommand{\co}{{\cal O}}
\newcommand{\smA}{\small {\A}}
\newcommand{\ttheta}{\tilde\theta}
\newcommand{\tn}{{\tilde\|}}
\newcommand{\rbar}{\overline{r}}
\newcommand{\oeps}{\overline{\varepsilon}}
\newcommand{\cgl}{\hbox{Lie\,}{\cal G}}
\newcommand{\Ker}{\hbox{Ker\,}}
\newcommand{\const}{\hbox{const.\,}}
\newcommand{\Sym}{\hbox{Sym\,}}
\newcommand{\tr}{\hbox{tr\,}}
\newcommand{\grad}{\hbox{grad\,}}
\newcommand{\ttd}{{\tt d}}
\newcommand{\ttdel}{{\tt \delta}}
\newcommand{\ns}{\nabla_*}
\newcommand{\csl}{{\cal SL}}
\newcommand{\kr}{\hbox{Ker}}
\newcommand{\beq}{\begin{equation}}
\newcommand{\eeq}{\end{equation}}
\newcommand{\pr}{\hbox{proj\,}}
\newcommand{\proj}{{\mathbb P}}
\newcommand{\tproj}{\tilde{\mathbb P}}
\newcommand{\projq}{{\mathbb Q}}
\newcommand{\tprojq}{\tilde{\mathbb Q}}
\newcommand{\oN}{\overline N}
\newcommand{\cN}{\cal N}
\newcommand{\cmet}{{{\hbox{${\mathcal Met}$}}}}
\newcommand{\met}{{\hbox{Met}}}
\newcommand{\bfga}{{\mbox{\boldmath$\overline\gamma$}}}
\newcommand{\bfOmega}{\mbox{\boldmath$\Omega$}}
\newcommand{\bfTh}{\mbox{\boldmath$\Theta$}}
\newcommand{\bfmw}{{\bf m}_{\mbox{\boldmath$\smo$}}}
\newcommand{\bfmu}{{\mbox{\boldmath$\mu$}}}
\newcommand{\bfulX}{{\mbox{\boldmath${X}$}}}
\newcommand{\bfultX}{\tilde{\mbox{\boldmath${X}$}}}
\newcommand{\bfmuw}{{\mbox{\boldmath$\mu$}}_{\mbox{\boldmath$\smo$}}}
\newcommand{\dmug}{d\mu_g}
\newcommand{\ltwo}{{L^{2}}}
\newcommand{\lfour}{{L^{4}}}
\newcommand{\hone}{H^{1}}
\newcommand{\honea}{H^{1}_{\smA}}
\newcommand{\er}{e^{-2\rho}}
\newcommand{\onetwo}{\frac{1}{2}}
\newcommand{\lra}{\longrightarrow}
\newcommand{\dv}{\hbox{div\,}}
\newcommand{\da}{D_{\mathbf A}}
\newcommand{\mbA}{\mathbb A}
\def\smo{{\mbox{\tiny$\omega$}}}
\protect\renewcommand{\theequation}{\thesection.\arabic{equation}}

\font\msym=msbm10
\def\Real{{\mathop{\hbox{\msym \char '122}}}}
\def\R{\Real}
\def\A{\mathbf A}
\def\Z{\mathbb Z}
\def\K{\mathbb K}
\def\J{\mathbb J}
\def\L{\mathbb L}
\def\D{\mathbb D}
\def\Mink{{\mathop{\hbox{\msym \char '115}}}}
\def\Integers{{\mathop{\hbox{\msym \char '132}}}}
\def\Complex{{\mathop{\hbox{\msym\char'103}}}}
\def\C{\Complex}
\font\smallmsym=msbm7
\def\smr{{\mathop{\hbox{\smallmsym \char '122}}}}
\def\ZHK{{\mathop{\hbox{\tiny ZHK}}}}

\maketitle
\thispagestyle{empty}
\vspace{-0.3in}
\begin{abstract}
We prove convergence of the gradient flow of the 
Ginzburg-Landau energy functional
on a Riemann surface in the self-dual Bogomolny case, 
in Coulomb gauge.
The proof is direct and makes use of 
the associated nonlinear first order differential operators
(the Bogomolny operators). One aim is to illustrate
that the Bogomolny structure, which is known to be of great
utility in the static elliptic case, can also be used effectively in 
evolution problems.
We also identify the minimizers
and minimum value of the energy when the 
Bogomolny bound is not achieved  (below the Bradlow limit).
\end{abstract}

\centerline{MSC classification: 58J35, 35Q56}

\section{{Introduction and statement of results}}
\label{secint}

We consider the self-dual Ginzburg-Landau, or abelian Higgs,
energy defined
over a two dimensional Riemann
surface $\Sigma$:
$$
\mcn({\mathbf A},\Phi)
= \frac{1}{2}\int_\Sigma 
D_ {\mathbf A}\Phi\wedge *D_{\mathbf A}\Phi+\;F *F + \frac{1}{4}\left(1
  - |\Phi|^2)\right)^2*1
$$
Here $\Phi$, the Higgs field, is a section of a complex line bundle
$L \to \Sigma$ with fixed metric $h$,  
and $D_ {\mathbf A}=(D_1, D_2)=\nabla -i\A $ is the covariant derivative operator 
defining an $S^1$ connection
on  $L$, with curvature 2-form $-iF$. We will assume that the degree
$N=deg\; L$ is a non-negative integer (but everything in this article has
an analogous version in the negative case).
It is known from \cite{brad88} that for $|\Sigma|>4\pi N$ the minimum
value of $\mcn$ is $\pi N$ and is attained on the set of self-dual vortices,
which are solutions of the Bogomolny equations \eqref{bogo}; these form a system
of first order partial differential equations, solutions of which are minimizers
of $\mcn$ {\em when they exist}. 
The presence of such a system of first order 
equations is a special feature of the functional $\mcn$ which is related
to self-duality (\cite{jt82,brad88}).
The main aim of this article is to show how it is possible to prove 
convergence of the corresponding gradient flow quite simply by making use of the Bogomolny (or self-dual) structure in \eqref{bogdec}-\eqref{bogo}; 
it is also proved that for $|\Sigma|<4\pi N$ the minimum value of 
$\mcn$ is $\pi N+\frac{l^2}{2|\Sigma|}$ where $l=2\pi N-\frac{|\Sigma|}{2}$,
and this is achieved by taking $\Phi\equiv 0$ and ${\mathbf A}$ to be a 
constant curvature connection.

A similar
proof of convergence for the $SU(2)$ Yang-Mills-Higgs 
functional on $\R^3$ was carried out in \cite{h93}; this is
closely related to the case $|\Sigma|>4\pi N$ in this paper. The
crucial structural features exploited there for Yang-Mills-Higgs also hold
for the Ginzburg-Landau flow in the case $|\Sigma|>4\pi N$: by \eqref{ebm}
the first order Bogomolny operators \eqref{bogo1}-\eqref{bogo2} converge
to zero, while \eqref{bom} then implies the convergence of  $(\A,\Phi)$ as $t\to +\infty$. 
The case $|\Sigma|<4\pi N$ is different:
indeed this latter inequality is an obstruction to the existence of solutions
to the Bogomolny equations (\cite{brad88}), 
and so the first stage of the argument just outlined
necessarily fails.  We give, however, an alternative set of quantities
in \eqref{ebm2}, whose large time behaviour can be analyzed, and which serve
as an effective alternative to the Bogomolny operators in this case.
Compared to the
adiabatic approximation method in \cite{ds97}, 
where convergence of the Ginzburg-Landau flow 
was proved for the case $\Sigma=
\R^2$, the proof here is certainly more efficient although it makes
less contact with the physics of vortices. It is also possible
to prove convergence by the very general method based on the
Lojasiewicz inequality in \cite{s83}, as has been done in \cite{ft}. 
In comparison with these previous results, the present method does provide 
the additional information that convergence is at an exponential rate.
Furthermore in view of the 
importance of the Bogomolny self-dual structure in the elliptic case,
it seems worthwhile to illustrate its utility
in the parabolic context, as was done also in \cite{ds09} for
vortex dynamics in a conservative (Hamiltonian) context.

Introducing conformal coordinates
$\{x^j\}_{j=1}^2$, in which the metric takes the form
$$
g=g_{jk}dx^j dx^k=e^{2\rho}
((dx^1)^2+(dx^2)^2)
$$ 
with associated area form $\dmug=e^{2\rho}dx^1\wedge dx^2$,
the energy  functional $\mcn$
is given by
\beq\la{fullenergy}
\mcn({\mathbf A},\Phi)
= \frac{1}{2}\int_\Sigma \Bigl(
|D_{\mathbf A}\Phi|^2e^{-2\rho} \;+\;B^2 + \frac{1}{4}\left(1
  - |\Phi|^2)\right)^2\Bigr) \ \dmug.
\ee
If we fix
a smooth connection $\nabla$ on $L$ then  $D_{\mathbf A}$ is determined uniquely by a
real 1-form 
$$ 
{\mathbf A}=A_1 dx^1 + A_2 dx^2 \in\Omega^1(\Sigma)
$$
where $A_i : \Sigma \to \Real$,  according to
\beq\la{fdc}
 D_\A :=\nabla - i{\mathbf A}.
\eeq
The curvature, or magnetic field $B$, of a connection
$D_{\mathbf A}=D_j dx^j$ is determined by $D_jD_k\Phi\,dx^j \wedge dx^k
=-iF\Phi=-iB\Phi d\mu_g$. Its integral
is a topological invariant of $L$:
\begin{equation}\label{topinv}
\int_\Sigma B \dmug = 2\pi N
\end{equation}
with $N= deg \;L$. 
We will choose $\nabla$ to have constant 
curvature $b$; this value is then fixed topologically as
$b=2\pi N/|\Sigma|$, where $|\Sigma|=\int_{\Sigma} d\mu_g$ is the
area of $\Sigma$. It follows that
$$B= b+*d{\mathbf A}\equiv b+e^{-2\rho}(\partial_1 A_2 - \partial_2 A_1).$$ 

\subsection{The Bogomolny equations}
A crucial property of the energy functional (\ref{fullenergy}) is that
(in appropriate functions spaces in which e.g.  integration by parts
is valid) it admits a {\it Bogomolny} decomposition into a sum of
squares of first order terms: 
\beq\la{bogdec}  
\mcn({\mathbf A},\Phi) =
\frac{1}{2}\int_\Sigma \Bigl( 4|\barpa \Phi|^2e^{-2\rho} \;+\;(B-
\frac{1}{2}\left(1 - |\Phi|^2)\right)^2\Bigr) \ \dmug \ + \ \pi N 
\eeq
where as above $N= deg L$ and
$$
\barpa \Phi = \frac{1}{2} (D_1 + i D_2)\Phi.
$$
If the following  first order equations, called the Bogomolny
equations, 
\begin{align}
\begin{split}
&\barpa\Phi=0,\\
B-&\frac{1}{2}(1-|\Phi|^2)=0\\
\end{split}
\la{bogo}
\end{align}
have solutions in a given class, they will automatically minimize
$\mcn$ within that class (this theory - and conditions under which
solutions exist - is developed for the plane and for surfaces in
\cite{jt82,nog,brad88}).     

Our main theorem concerns the large time behaviour of global weak
(Sobolev) solutions to  
\ref{eom} in Coulomb gauge (the existence and uniqueness of which is shown in section \ref{secexist}).   We show that as $t\to\infty$ this solution 
converges exponentially fast to a {\it minimizer} of  \eqref{fullenergy} 
a solution of the 
{\it Bogomolny equations} when $|\Sigma| >4\pi N$ .   We make the assumption of
small (close to equilibrium) initial energy. 

We will make use of the Bogomolny equations to
derive the asymptotic convergence of weak solutions to the gradient
flow equations in section \ref{secasy}.  In fact, the convergence of
such solutions to minimizers of the energy provides an independent
proof of the existence of minimizers and solutions to the Bogomolny
equations in the case when $|\Sigma|>4\pi N$.

\subsection{The time-dependent  equations and statement of the main theorem}
Associated with the energy functional \eqref{fullenergy}  are the  Ginzburg-Landau
gradient flow equations in the variables 
$(\mbA, \Phi)\equiv (A_0, \A,\Phi)\equiv (A_0, A_1, A_2 , \Phi)$ on
$[0,\infty)\times \Sigma$:
$$
(\frac{d}{dt}, D_0)(\A, \Phi) = - \mcn(\A,\Phi)
$$
and equivalently,
\begin{align} 
\begin{split}
  & {\frac{\partial A_1}{\partial t}}-
{\frac{\partial A_0}{\partial x^1}}=
- \frac{\partial B}{\partial x^2}+
\langle i\Phi ,{D}_1\Phi\rangle,\\
  & {\frac{\partial A_2}{\partial t}}-
{\frac{\partial A_0}{\partial x^2}}=
+ \frac{\partial B}{\partial x^1}+ 
\langle i\Phi ,{D}_2\Phi\rangle,\\ 
  & D_0\Phi -\Delta_{\mathbf A}\Phi
        -\frac{1}{2}(1-|\Phi|^2)\Phi=0\,,
\end{split}\label{eom}
\end{align}
where the covariant Laplacian $\Delta_{\mathbf A}$ is defined as
$$
\Delta_{\mathbf A} = e^{-2\rho} ( D_1^2 +D_2^2 )\,.  
$$
This is the fully gauge invariant version of the Ginzburg-Landau
gradient flow. The Higgs field
$\Phi$, now depending on time $t\geq 0$ , is a section of a complex line bundle
$\L\equiv [0,\infty)\times L \to [0,\infty)\times \Sigma$   
on which there is a connection 
$$
{\mathbb A} = (A_0, A_1,A_2) = (A_0, {\mathbf A})
$$ 
with associated covariant derivative 
$\D = (D_0, D_1, D_2) = (D_0, D_\A)$,
where $D_0 = \pt - i A_0 $ is the time component and 
${D_\A}=(D_1, D_2)$ is the spatial component as in \eqref{fdc}
with  $A_i$ now depending on time $t\geq 0$ as well as $x\in\Sigma$.   
The dependent variables are thus the Higgs field $\Phi$ and a real 1-form
$A_0dt+A_j dx^j$  on $[0,\infty)\times\Sigma$. The system is invariant
under the infinite dimensional group of gauge transformations:
if $\chi(t,x)$ is a smooth function then $(A_0,A_1,A_2,\Phi)$ is a smooth
solution of \eqref{eom} 
if and only if 
$(A_0+\partial_t\chi,A_1+\partial_1\chi,A_2+\partial_2\chi, e^{i\chi}\Phi)$
is a smooth solution of \eqref{eom}. This means it is possible to place
further conditions on the solution, and this is necessary to obtain
uniqueness. We impose the  {\it Coulomb  gauge condition} :
\beq
\dv {\mathbf A}=e^{-2\rho}(\partial_1 A_1+\partial_2 A_2)=0 \ \ ;
\la{g1}\eeq 
a solution, or set of initial data, is said to be in {\it Coulomb gauge} 
when this condition holds for all relevant $t$. In addition, we will 
require that 
\beq\int_{\Sigma} A_0 d\mu_g=0
\la{g2}\eeq at all relevant $t$: this  may be achieved by applying to 
the solution 
$(A_0, A_1,A_2, \Phi)$ at each $(t,x)$ the gauge transformation   
$\chi (t) = -t|\Sigma|^{-1}\int_\Sigma A_0 \dmug$.
\\[1ex]

{\it Initial data and Coulomb gauge}: We specify
$({\mathbf A},\Phi)= (A_1, A_2,\Phi)$ at $t=0$
which in Coulomb gauge also determines $A_0$ initially.  To see this,
we eliminate $A_0$ from the full (gauge invariant) equations \eqref{eom} 
by taking the divergence in the equation for ${\mathbf A}$  under the 
condition \eqref{g1}.  This leads to the equation for $A_0$ in Coulomb gauge,
\beq 
-\Delta A_0=\dv \langle i\Phi,D_{\mathbf A}\Phi\rangle \la{a0eq} \  .
\eeq 

\begin{lemma}\la{a0est}
For each fixed $t\geq 0$ there  is a smooth function  $({\mathbf A},\Phi) \mapsto \alpha({\mathbf A},\Phi)$ 
from  $H^1(L\oplus\Omega^1)\to H^{1,q}$, for any $1<q<2$, 
which is
the unique solution of \eqref{a0eq} satisfying \eqref{g2}.
\end{lemma}
\proof
Because $\Phi\in H^1$ then for all $p<\infty$,  $\Phi\in L^p$, 
$D_{\mathbf A}\Phi\in L^2$ and so $\dv \langle i\Phi,D_{\mathbf A}\Phi\rangle\in
H^{-1,q}$ for all $q<2$. It follows that
\eqref{a0eq} has a unique solution $A_0 =  \alpha({\mathbf
  A},\Phi)$ verifying \eqref{g2} as described in the lemma. 
\qed
\begin{remark}
Furthermore, if $\Phi$ satisfies \eqref{eom} then for $t\  a.e.$
$
-\Delta A_0=\dv \langle i\Phi,D_{\mathbf A}\Phi\rangle
=\langle i\Phi,D_0\Phi\rangle\,,
$
and so by standard elliptic theory, 
$$
\|A_0\|_{H^{2,p}} \leq c_s \|\langle i\Phi , D_0\Phi\rangle\|_{L^{p}}\,,
$$
where finiteness of the right hand side will be justified in 
section \ref{34}.
By Sobolev imbeddings, here applied for $1<p<2$ at each time
$t\geq 0$,
\beq\label{estao}
\|A_0\|_{H^{2,p}}
\leq C_1\|\Phi\|_{L^{\frac{2p}{2-p}}}\|D_0\Phi\|_{L^2}\quad\hbox{ and so}
\quad\|dA_0\|_{L^2}\leq C_2\|\Phi\|_{\honea}\|D_0\Phi\|_{L^2}.
\eeq

\end{remark}
As a particular consequence,  $A_0$ is initially determined through the 
initial data for $(\A, \Phi)$.  In addition, the estimates \eqref{estao}  
will be used below to derive  estimates for $D_0\Phi$.

Explicitly applying the gauge conditions \eqref{g1} and \eqref{g2}
and substituting $A_0=\alpha({\mathbf A}, \Phi)$,
the system \eqref{eom} becomes the following nonlocal parabolic system:
\begin{align} 
\begin{split}
\dot {\mathbf A} & =\Delta {\mathbf A}+d\alpha+\langle i\Phi,D_{\mathbf A}\Phi\rangle\\
\dot\Phi & =\Delta\Phi-2ie^{-2\rho}\sum_jA_j\partial_j\Phi
-e^{-2\rho}|{\mathbf A}|^2\Phi + i \alpha\Phi+\frac{1}{2}(1-|\Phi|^2)\Phi \ 
\end{split}
\la{eomc}
\end{align}
with initial data specified for  $(\A, \Phi)$.
We discuss existence of this system in section \ref{secexist} and we will need the following spaces.

{\it Function spaces}:
We will work with  spaces of functions $H^{s,p}(\Sigma)$ which may be defined
either 
(i) by using a partition of unity 
to reduce to the Euclidean case (the standard spaces
$H^{s,p}(\R^2)$), or
(ii) by appealing to the spectral theory for $-\Delta$
to define a functional calculus and then
introducing the norm as $$\|f\|_{H^{s,p}}=\|(1-\Delta)^{\frac{s}{2}}f\|_{L^p}.$$
Let $f_j$ be an orthonormal basis of eigenfunctions,
$-\Delta f_j=\lambda_j f_j$ with $0\leq \lambda_1\leq\lambda_2\leq\dots$,
then for $p=2$ an equivalent norm on $H^s=H^{s,2}$ is
\beq
\|f\|_{H^s}^2=\sum_j (1+\lambda_j)^s|c_j|^2,\quad\hbox{where}\;f=\sum c_j f_j.
\label{defn}
\eeq

These definitions  extend to sections of vector bundles $L$, in which
case we write $H^s(L)$ for the corresponding spaces defined with
respect to the connection  $\nabla$ and its corresponding Laplacian
$e^{-2\rho} ( \nabla_1^2 +\nabla_2^2 )$; for $s=1$ an equivalent norm is
\beq
\la{dhonea}
\|\Phi\|^{2}_{H^1(L)} = 
\int_\Sigma\left( |\Phi|^2 + |\nabla \Phi|^2 \right)d\mu_g . 
\eeq
In the above integral  the inner products are the standard ones
induced from $h$ and $g$. 
The definitions also extend  to 1-forms 
using the Hodge Laplacian $\Delta=-(d\delta+\delta d)$,
(with the sign chosen so that it is non-positive). 
In conformal coordinates:
\begin{align}
(\Delta {\mathbf A})_i & =\partial_i(e^{-2\rho}(\partial_1 A_1+\partial_2 A_2))
-\epsilon_{ij}\partial_j(e^{-2\rho}(\partial_1 A_2-\partial_2 A_1))\cr
& = -\epsilon_{ij}\partial_j B
\notag\end{align}
when the Coulomb gauge \eqref{g1} holds. The corresponding norms are written
$H^s(\Omega^1)$ where $\Omega^1$ is  the space of one-forms on $\Sigma$ with 
coefficients in $H^s$.\\
\\
{\it The main theorem}\\
For $H^1$ data we  derive existence of a solution also in $H^1$ in the following
section in the Coulomb gauge.  The main result concerns the time asymptotic behaviour of these solutions which are shown to be {\it exponentially} 
converging to  an energy minimizer.   The minimum  energy
depends on whether the size of $\Sigma$ is large enough to support the
existence of vortices (\cite{brad88,nog}):
$$ \mcn_{min} ={ \left\{
   \begin{split}
     \pi N \phantom{ + \frac{(4\pi N-{|\Sigma|})^2}{8|\Sigma|}}& \mbox{ if } \Sigma > 4\pi N \\
     \pi N + \frac{(4\pi N-{|\Sigma|})^2}{8|\Sigma|} &\mbox{ if } \Sigma \leq 4\pi N \ .
   \end{split}
\right\} 
}
$$
Let $$
\mcn_0 = \mcn ( \A(0),\Phi(0))
$$
be the initial energy.
The following is the main theorem 

\begin{theorem}[Main Theorem]\label{main}
Given finite energy initial data $(\A(0),\Phi(0))\in
H^1(\Omega^1 \oplus L)$ satisfying the Coulomb gauge condition
\eqref{g1} and such that $\epsilon_0$ with 
$$0\ \leq \ \mcn_0 - \mcn_{min} \ =\ \epsilon_0 $$ is sufficiently small  there exists a unique global solution
$(A_0,\A,\Phi)$  with $A_0 \in C ([0,\infty)); H^{1,p}(\Sigma))$
for any $p<2$, and $(\A , \Phi)\in C\big([0,\infty); H^1(\Omega^1
\oplus L)\big)$ of \eqref{eom},  satisfying \eqref{g1} and \eqref{g2} with the
following large time behaviour:
\begin{itemize}
\item[] if $ |\Sigma| \ > \ 4\pi N$ then $({\mathbf A},\Phi)
  \stackrel{t\to \infty}{\longrightarrow} ({\mathbf
    A}_\infty,\Phi_\infty)$ in $L^2(\Omega^1\oplus L)$ where
  $({\mathbf A}_\infty,\Phi_\infty)$ is a solution of the minimum energy static
  (Bogomolny) equations and where $\Phi_\infty$ has precisely $N$
  zeros (vortices) on $\Sigma$.  Also  $A_0(t)\to 0$ strongly in $L^q$
  for all  $q<\infty$.   
The minimum value of the energy in $H^{1}(L \oplus \Omega^1)$ is attained and  is $\pi N$.
\item[]
if $ |\Sigma| \ < \ 4\pi N$ then $({\mathbf A},\Phi)  \stackrel{L^2}{\to} ({\mathbf A}_\infty,0)$ 
where
$D=\nabla-i{\mathbf A}_\infty $ is a constant curvature connection at which the minimum value of the energy in $H^1 (L\oplus \Omega^1)$ is achieved and equals
$$
\mcn({\mathbf A}_\infty ,0)\,=\,\min_{H^1(L\oplus\Omega^1)}\mcn= \pi N+\frac{(4\pi N-{|\Sigma|})^2}{8|\Sigma|}\,.
$$
\end{itemize}
Convergence is at an exponential rate, namely,
$\exists \ c, \delta>0$ depending only  the initial data such that
\beq\la{converg}
\|({\mathbf A}, \Phi) - ({\mathbf A}_\infty, \Phi_\infty)\|_{L^2} \leq c\ e^{-\delta t},\quad \mbox{respectively,}\quad
\|({\mathbf A}, \Phi) - ({\mathbf A}_\infty, 0)\|_{L^2} \leq c\ e^{-\delta t}  .
\eeq
\end{theorem}

Global existence for \eqref{eom} is established in the following section, making use of  standard (semigroup) techniques (theorem \ref{exist}).  The proof of asymptotic convergence, which is the main content of theorem \ref{main},
is then given in section \ref{secasy}.

\begin{remark}
Convergence in stronger $H^s$ norms for $s<1$
can be deduced by interpolation between $L^2$ and $H^1$ from \eqref{converg}, using convergence in $L^2$ (``low'' norm) and boundedness in $H^1$ (or ``high'' norm) 
guaranteed  the energy non-increase.
Furthermore, 
when boundedness can be proved in higher norms than $H^1$ then
interpolation would lead to convergence in $H^s$ for higher $s$ by standard parabolic theory.
\end{remark}
\begin{remark}
The case 
$$
|\Sigma|=4\pi N
$$
is degenerate (as becomes clear from reading the proof) and it is 
conceivable that exponential convergence may not occur, even if the initial 
energy is close to the minimum energy.  
\end{remark}
\begin{remark}
The condition that the  initial energy be  
close to the minimum energy can be relaxed in situations when
there are no non-minimal critical points (\cite[chapter 3]{jt82}).
The energy is 
strictly decreasing (as it is a gradient flow) and so it must
decrease to the energy of a critical point. 
If non-minimal critical points do not exist then the energy must 
decrease to its minimum value,
in which case it approaches the minimum energy,  and then the
conclusion of the theorem holds 
regarding exponential convergence.  This is not discussed here further. 
\end{remark}

\section{Existence theorem}\label{secexist}

There are different ways to obtain a global existence theorem 
for these equations, for example using  maximum principles 
(as for the case of the same equations
on $\R^2$ in \cite{ds97}) or energy methods as done here, 
using semi-group techniques explained in  e.g. \cite{taylor} and elsewhere.   
To use the energy norm method, we first
reduce our system in Coulomb gauge to a system for $(\A,\Phi)$ given in \eqref{eomc}; 
then we  show that  for $H^1$ data $({\mathbf A}(0),\Phi (0))$  there is a unique global weak solution of 
those equations,  also in $H^1$, varying continuously in time.  

\begin{theorem}[Existence in Coulomb gauge]
\la{exist}
Given initial data $({\mathbf A}(0),\Phi(0))\in H^1(\Omega^1 \oplus
L)$
(which implies finite initial energy
$\mcn(\A(0),\Phi(0)) =\mcn_0 < \infty$),
there exists a unique mild solution of \eqref{eomc} with
regularity $({\mathbf A}, \Phi)\in C([0,\infty);H^1(\Omega^1 \oplus L)$,
which  satisfies
$$\mcn(({\mathbf A}, \Phi)(t))\leq\mcn(({\mathbf A}, \Phi)(0))$$ 
for
all $t\geq 0.$
For smooth initial data the solution is smooth, and for
$H^1(\Omega^1 \oplus L)$ initial data the solution is the limit
in $C_{loc}([0,\infty);H^1(\Omega^1 \oplus L)$ of smooth solutions.
\end{theorem}
\proof
Writing $\cu=({\mathbf A}, \Phi)$ and $\cl\cu=(\Delta {\mathbf A}, \Delta\Phi )$
the system \eqref{eomc} is of the form $\dot\cu= \cl\cu+\mcF(\cu)$,
which can be treated as
an abstract evolution equation in the space $X=H^1(\Omega^1 \oplus L )$.
We introduce also the auxiliary space $Y=H^{-\frac{1}{2}}
(\Omega^1\oplus L)$, and make use of the following facts:
\begin{itemize}
\item
$e^{t\cl}$ is a strongly continuous semi-group of  contractions 
on $X$ for $t\geq 0$,
\item
$e^{t\cl}:Y\to X$ for $t>0$, with $\|e^{tL}\|_{Y\to X}\leq C t^{-3/4}$ for $0<t\leq 1$,
\item 
$\mcF:X\to Y$ is a smooth function satisfying  $\|\mcF(\cu_1)-\mcF(\cu_2)\|_{Y}
\leq K(R)\|\cu_1-\cu_2\|_X$ for $\|\cu_j\|_X\leq R$.
\end{itemize}
The first of these is a standard property of the heat equation. The
second can be derived using the spectral representation introduced above:
$$e^{t\Delta}\sum c_j f_j=\sum e^{-t\lambda_j}c_j f_j$$
which implies that
$$
\|e^{t\Delta}f\|_{H^s}^2=\sum e^{-2t\lambda_j}
(1+\lambda_j)^s|c_j|^2\leq 
C_{s,r}(t)^2\|f\|_{H^r}^2
$$
with $C_{s,r}(t)=\sup_{x\geq 0}\frac{(1+x)^{\frac{s}{2}}e^{-tx}}{(1+x)^{\frac{r}{2}}}\sim
c t^{-\frac{s-r}{2}}$ as $t\to 0$.
The third assertion follows by examining the various terms which constitute
$\mcF$ and applying appropriate embeddings for the $H^{s}$ spaces, in particular
$H^{\frac{1}{2}}\subset L^4$ and $L^{\frac{4}{3}}\subset H^{-\frac{1}{2}}$
(continuous embeddings). For example, the current term 
$\langle i\Phi, \nabla\Phi\rangle$ in the $\Phi$ equation
arises as a continuous bilinear map $H^1\times L^2\to L^{\frac{4}{3}}\subset 
H^{-\frac{1}{2}}$ bounded by
$|\langle i\Phi, \nabla\Phi\rangle|_{L^{\frac{4}{3}}}
\leq c\|\Phi\|_{L^4}\|\nabla\Phi\|_{L^2}$. The other terms in $\mcF$ are treated
similarly (using lemma \ref{a0est} to handle the $\alpha$ terms).

These three properties imply (see \cite[\S15.1]{taylor}
that the integral operator
$\cu\to\int_0^t e^{(t-s)\cl}\mcF(\cu(s))ds$ is a contraction
on $C([0,T];X)$ for sufficiently small $T$, and in fact $T$
can be taken to be a positive non-increasing function of
$\|\cu(0)\|_X$. This implies that there is a unique local
solution $\cu\in C([0,T];X)$, and also that, given two sets
of initial data $\cu_1(0),\cu_2(0)$ in $X$, the well-posedness estimate 
$$
\max_{0\leq t\leq T}\|\cu_1(t)-\cu_2(t)\|_X\leq c\|\cu_1(0)-\cu_2(0)\|_X
$$
holds for some $c>0$, for $T$ sufficiently small (again depending
only on $\|\cu_j(0)\|_X$.)  From this it can be deduced that for
$H^2$ initial data the solution remains in $H^2$, and in fact for
$H^s$ initial data the solution remains in $H^s$ for all $s\geq 2$.
Thus the $H^1$ solutions can be approximated by regular solutions,
and obey the energy non-increase :
$\mcn(({\mathbf A}, \Phi)(t))\leq\mcn(({\mathbf A}, \Phi)(s))$ for $t\geq s$. As
a consequence the $H^1$ norm is globally bounded
and there 
exists a unique global solution $\cu\in C([0,\infty);X)$ of the
corresponding integral equation
$$
\cu(t)=e^{t\cl}\cu(0)+\int_0^t e^{(t-s)\cl}\mcF(\cu(s))ds
$$
with the approximation property asserted.
\qed

\section{Asymptotic behaviour and proof  of the main theorem}
\label{secasy}
We continue  to use the  Coulomb gauge in which $A_0$ is eliminated as
an independent variable by lemma \ref{a0est}.
We consider the cases $|\Sigma|> 4\pi N$ and $|\Sigma|<4\pi N$ separately.
In each case we introduce auxiliary variables which are specially tailored
to reveal the asymptotic convergence of $(\A, \Phi)$ 
to a limit $(\A_\infty,\Phi_\infty)$ which is characterized differently 
depending on whether the surface area $|\Sigma|$ is bigger or smaller than
$4\pi N$. In the following proofs we shall make use of the norm
\beq
\la{honea}
\|\eta\|^{2}_{\honea} = 
\int_\Sigma\left( |\eta|^2 + |D_{\A}\eta|^2 \right)d\mu_g 
\eeq
defined with respect to the time dependent connection $D_{\A}=\nabla-i\A$.
\subsection{ Proof in the case $|\Sigma|> 4\pi N$}
\label{34}

It is useful to introduce   the Bogomolny variables
$v=v(\A,\Phi)$ and $\eta=\eta(\A,\Phi)$, defined as 
\begin{align}
\eta &= \barpa\Phi\la{bogo1} \\
v&= B-\frac{1}{2}(1-|\Phi|^2)\ \la{bogo2}
\end{align}
where, as above, $\barpa = \frac{1}{2} (D_1 + i
D_2)$.
In terms of these variables the energy decomposes:
$$
\mcn =
\frac{1}{2}\int_\Sigma\,( 4|\eta|^2e^{-2\rho}
+v^2)\, \dmug \ + \ \pi N \,,
$$
as was seen above in  \eqref{bogdec}.  
The evolution of these variables is  according to
\begin{align}\la{ebm}
\begin{split}
&D_0\eta-4\barpa(e^{-2\rho}\pbfa\eta)+|\Phi|^2\eta=-v\eta,\\
&(\pt-\Delta+|\Phi|^2)v=-4|\eta|^2.
\end{split}
\end{align}
These equations can be obtained from  \eqref{eom} by first applying the operator
$\barpa$ to the equation for $\Phi$,
using the identities 
(in conformal co-ordinates)
$$
4\pa\barpa = D_{1}^{2} + D_{2}^{2} + e^{2\rho} B\,,
\qquad\qquad
\Delta_{\mathbf A}=e^{-2\rho}(D_1^2+D_2^2)
$$
and
\begin{align*}
&\barpa\Delta_{\mathbf A} -4\barpa(e^{-2\rho}\pa)\barpa
= - (\bar\partial B) - B\barpa 
\intertext{and}
&[\barpa ,  D_0] = \frac{i}{2}(\dot{A}_1 + i\dot{A}_2) - i \bar\partial A_0
\end{align*}
(where $\bar\partial A_0 = (\bar\partial A_0) + A_0 \bar\partial $) 
and finally using the equation for ${\mathbf A}$ in the last
commutator. (To derive the final commutator use also  that the 
background connection $\nabla$ is fixed independent of $t$, i.e.,
$[\pt, \nabla] = 0$.)

We will show (using lemma \ref{lb}) that 
$$
(\eta, v)\longrightarrow (0,0) \mbox{ as } t\to \infty
$$
(in specified spaces) and then deduce via lemma \ref{aphidot}
the convergence of
$({\mathbf A}, \Phi)$ as stated in the main theorem.   
For this consider two energy-type quadratic forms,
respectively associated to  each of the equations in \ref{ebm} (at each time
$t>0$) defined as  
\begin{align}
&Q_{\A,\Phi}(\eta)\;=\;\int_{\Sigma}\Bigl(
4e^{-4\rho}|\pbfa\eta|^2+e^{-2\rho}|\Phi|^2|\eta|^2
\Bigr)d\mu_g \ ,\la{qfeta}\\
&Q_\Phi(v)\;\;\;\,\,=\;\int_\Sigma\Bigl(e^{-2\rho}(|D_1 v|^2+|D_2 v|^2)
+|\Phi|^2v^2\Bigr)d\mu_g \la{qfv} .
\end{align}
\begin{remark}
It is shown in   \cite{ds09}[lemma 3.2.2] that these quadratic forms are  
coercive
in $\honea$ and  $H^1$ respectively, provided that $\|\Phi\|_{L^2}\geq m >0$:
to be precise, under this assumption on $\Phi$, there  exists
$\gamma=\gamma(m,M)>0$ (where recall that $\mcn \leq M <\infty$ and $M=\pi N +1$),  such that 
\beq\la{lbs}
Q_{\A,\Phi}(\eta)\geq\gamma\|\eta\|_{\honea}^2\quad\hbox{and }
\quad
Q_\Phi(v)\geq\gamma\|v\|_{H^1}^2\,.
\eeq   
That the assumption on $\|\Phi\|_{L^2}$
is a valid one  can be shown as follows:
by the Bogomolny decomposition of the functional $\mcn$
\eqref{bogdec} and the topological invariant \eqref{topinv},
$$
\int v = 2\pi N -\frac{|\Sigma|}{2} +\frac{\|\Phi\|^2_{L^2}}{2}
$$
from which it follows that, 
\begin{align}
\begin{split}\label{vpsn}
\|\Phi\|_{L^2}^{2} 
  & = \ |\Sigma| - 4\pi N  +2\int v \ \\
  & \geq \ |\Sigma| - 4\pi N  - 2
(2\epsilon_0)^{\frac{1}{2}}|\Sigma|^{\frac{1}{2}} \\
  & > 0 
\end{split}
\end{align}
for  $\epsilon_0$ sufficiently  small.
\end{remark}
\medskip
In terms of the quadratic forms $Q_{\A,\Phi}, Q_{\Phi}$ the variables $\eta, v$ satisfy the following inequalities:

\begin{lemma}\la{lb}
Let $(\A,\Phi)$ be as in theorem \ref{exist} and define $(\eta,v)$
as in \eqref{bogo1}-\eqref{bogo2}. Then
for any time $T>0$, the following identities hold
\begin{align}
\la{enid}
&e^{2\delta T}\|\eta(T)\|_{L^2}^2+2\int_0^Te^{2\delta t}\Bigl(
Q_{\A,\Phi}(\eta)-\delta\|\eta\|_{L^2}^2
\Bigr)\,dt\;\leq \;\|\eta(0)\|_{L^2}^2-2\int_0^T\int_\Sigma e^{2\delta t}v|\eta|^2 
e^{-2\rho}d\mu_g dt,\\
\la{enid1}
&e^{2\delta T}\|v(T)\|_{L^2}^2+2\int_0^Te^{2\delta t}\Bigl(
Q_{\Phi}(v)-\delta\|v\|_{L^2}^2
\Bigr)\,dt\;\;\;\;\leq\;\|v(0)\|_{L^2}^2
-8\int_0^T\int_\Sigma e^{2\delta t}v|\eta|^2e^{-2\rho}d\mu_g dt.
\end{align}
\end{lemma}

\proof
For smooth solutions multiply \ref{ebm}  respectively by $\eta $ and
$v $ and integrate with respect to $e^{2\delta t}e^{-2\rho} \dmug dt$ 
over  $[0,T]\times\Sigma$.
For more general finite energy solutions
use the approximation property in theorem \ref{exist}. 
\qed
\\[2ex]

  The above lemma
implies the exponential decay of $\eta, v$ under the conditions on $\Phi$ for 
co-ercivity and for initial energy sufficiently close to its minimum.

\begin{corl}[Exponential decay]
For all $\delta\in (0,\frac{\gamma}{2})$ and  $0< \epsilon_0 $ with $ \sqrt{2\epsilon_0} \leq \frac{\delta}{8 C_s^2}$,  if 
$\mcn(A,\Phi) -\pi N \leq \epsilon_0$,  then 
\beq\la{int}
\sup_{t>0}e^{2\delta t}(\|\eta(t)\|_{L^2}^2 +\|v(t)\|_{L^2}^2) +\delta \int_0^\infty
e^{2\delta t}(\|\eta\|_{\honea}^{2} +\|v\|_{\hone}^2)\,dt 
\leq\ 2\bigl(\|\eta(0)\|_{L^2}^2 
+\|v(0)\|_{L^2}^{2}\bigr)\,.
\eeq
\end{corl}
\proof
Recombining the terms in \eqref{enid} we have the following bounds
\begin{align}
  & 2\int_0^Te^{2\delta t}\Bigl(
 Q_{\A,\Phi}(\eta)-\delta\|\eta\|_{L^2}^2
 \Bigr)\,dt\  + 2\int_0^T\int_\Sigma e^{2\delta t}v|\eta|^2 
 e^{-2\rho}d\mu_g dt \notag \\
  &   \phantom{\int_0^Te^{2\delta t}\Bigl(Q_{\mbA,\Phi}
            (\eta)-\delta\|\eta\|_{L^2}^2 \Bigr)\,dt\  
           }
     \geq 2\int_0^T e^{2\delta t}(\gamma-\delta)\|\eta\|_{\honea}^{2}
   \ dt \ -\ 2\sup_{t}\|v\|_{L^2}  \int_0^T e^{2\delta
     t}\|\eta\|_{L^4}^{2} \ dt \notag \\
  &   \phantom{\int_0^Te^{2\delta t}\Bigl(
            Q_{\mbA,\Phi}(\eta)-\delta\|\eta\|_{L^2}^2 \Bigr)\,dt\  
           }
      \geq \int_0^T e^{2\delta t}(2\gamma-2\delta\;-\; 4C_s^2\sup_t
    \|v\|_{L^2}  )\|\eta\|_{\honea}^{2} 
   \ dt \la{lineqb1}
\end{align}
as there exists $C_s>0$ such that $\|\eta\|_{L^4}
\leq C_s\|\eta\|_{\honea}$ by the Sobolev inequality.
(Covariant Sobolev spaces and extensions of standard inequalities in
these were discussed in the appendix in \cite{ds09}).  By the energy non-increase,
\be\label{lastineq}
\sup_{[0,T]}\|v\|_{L^2}
\leq \sqrt{2({\mcn} - \pi N  )} \ \leq \
\sqrt{2({\mcn}_0 - \pi N  )} \ = \
\sqrt{2\epsilon_0}
\ee
and so for
$\sqrt{2\epsilon_0}
\leq
\frac{2\gamma-3\delta}{4C_s^2} $ (which is implied by the assumption on $\epsilon_0$),
the final term on the right hand side of \eqref{lineqb1} is bounded below as
$$
\geq \delta \int_0^T e^{2\delta t}\|\eta\|_{\honea}^{2}\,dt\,.
$$
Hence from \eqref{enid} 
\beq\la{int2}
\sup_{t\in [0,T]}e^{2\delta t}\|\eta(t)\|_{L^2}^2  \ 
+\ \delta \int_0^T e^{2\delta t}\|\eta\|_{\honea}^{2}\,dt 
\leq\ \|\eta(0)\|_{L^2}^2 
\,.
\eeq
Using this and \eqref{lastineq} (together with Holder's inequality and the same 
Sobolev imbedding as above)
we can bound the final term on the right hand side of \eqref{enid1}:
$$
8\int_0^T\!\int_\Sigma e^{2\delta t}v|\eta|^2 
 e^{-2\rho}d\mu_g dt
\leq
8\int_0^T\!\int_\Sigma e^{2\delta t} \|v\|_{L^2}\;  C_s^2\|\eta\|_{\honea}^{2} \ dt
\leq
\frac{8C_s^2\sqrt{2\epsilon_0}}{\delta}
\|\eta(0)\|_{L^2}^2
\leq\|\eta(0)\|_{L^2}^2
$$
for  $\sqrt{2\epsilon_0} \leq \frac{\delta}{8 C_s^2}$.
Using this and   \eqref{int2} in \eqref{enid2} we have 
\beq\label{int3}
\delta \int_0^T e^{2\delta t}\|v\|_{\hone}^{2}\,dt\,
\leq 2\int_0^Te^{2\delta t}\Bigl( Q_{\Phi}(v)-\delta\|v\|_{L^2}^2
 \Bigr)\,dt\  
\leq\ \bigl(\|\eta(0)\|_{L^2}^2+\|v(0)\|_{L^2}^2\bigr)\,
\eeq
which proves the corollary.
\qed

In terms of $v=B-\frac{1}{2}(1-|\Phi|^2)$ and
$\eta = \barpa\Phi=\frac{1}{2}(D_1+iD_2)\Phi$, the
equations \eqref{eom} read (in gauge invariant form):
\begin{align}
\begin{split}
  & {\frac{\partial A_1}{\partial t}}-
{\frac{\partial A_0}{\partial x^1}}=
- 2
\langle \Phi ,i\eta\rangle
-{\frac{\partial v}{\partial x^2}},\\
  & {\frac{\partial A_2}{\partial t}}-
{\frac{\partial A_0}{\partial x^2}}=
-2\langle \Phi ,\eta\rangle
+\frac{\partial v}{\partial x^1}
,\\ 
&D_0\Phi=4e^{-2\rho}\pa\eta-v\Phi.\\
\end{split}
\label{bom}
\end{align}

We will now show that the equations \eqref{bom} together with the
estimates \eqref{int} and also \eqref{estao} for $A_0$ 
imply  the convergence as $t\to \infty$
in $L^2$ of $t\mapsto (\A,\Phi)(t)$ at an exponential
rate.
\begin{lemma}\la{aphidot}
In terms of $\eta, v$ above,
\begin{equation}\la{dot}
\|\dot {\mathbf A}\|_{L^2}+\|\dot\Phi\|_{L^2}\leq c(
\|\eta\|_{\honea}+\|v\|_{H^1})
\end{equation}
which together with the estimate \eqref{int} completes the proof of the 
theorem in the case $|\Sigma|> 4\pi N$.
\end{lemma}
\proof
Unless specified otherwise the generic constants $c,c'$ depend only on $\Sigma$.

Firstly, 
$$
\int_\Sigma (1-|\Phi|^2)^2  = |\Sigma| - 2\|\Phi\|^{2}_{L^2}  +
\|\Phi\|^{4}_{L^4}\geq |\Sigma| -
2\|\Phi\|_{L^4}^{2}|\Sigma|^{\frac{1}{2}} +\|\Phi\|^{4}_{L^4}\geq
\frac{1}{2}\|\Phi\|^{4}_{L^4} -|\Sigma|. 
$$  
(as  $\|\Phi\|^{2}_{L^2}\leq \|\Phi\|_{L^4}^{2} |\Sigma|^{\frac{1}{2}} \leq \frac{1}{4}  
\|\Phi\|_{L^4}^{4} + |\Sigma|$).  
Therefore,
$$
\mcn ({\mathbf A},\Phi) \ \geq \ c\big(\|D_{\mathbf A} \Phi\|_{L^2}^{2} +
\|\Phi\|_{L^4}^{4}\bigr)-c'  \ \geq \  c\|\Phi\|_{\honea}^{2}-c'
$$ 
and using  the energy non-increase (shown at the end of the proof 
of the existence theorem \ref{exist}) this  implies
\begin{equation}\la{phil4}
\|\Phi\|_{L^4}^{4} \leq c(\mcn_0 , \Sigma)
\end{equation} 
where $\mcn_0$ is the initial energy.  

The last equation in \eqref{bom} implies 
\begin{align*}
\|D_0\Phi\|_{L^2}^{2} 
  & \leq \;c\big(\|\eta\|^{2}_{\honea} + \|v\Phi\|_{L^2}^{2}\big) \\
  & \leq \;c\big(\|\eta\|^{2}_{\honea} +
  \|\Phi\|_{L^4}^{2}\|v\|_{L^4}^{2}\big) \\
  & \leq c(\mcn_0, \Sigma) \big(\|\eta\|^{2}_{\honea} +
    \|v\|_{H^1}^{2}\big)
\end{align*}
Therefore by \eqref{int},
\be\la{dophi}
\int_{0}^{\infty} \|D_0\Phi\|^{2}_{L^2} e^{2\delta t}\ dt \leq \infty
\ee
which implies the exponential decay of $t\mapsto D_0 \Phi$.  
However, as seen in the previous section,   the equation for $\Phi$
implies in the Coulomb gauge
$$
-\Delta A_0 =\langle i\Phi,\ D_0\Phi\rangle 
$$
and then the Calderon-Zygmund inequality, given the condition $\int A_0 \dmug = 0$, implies
$$
\|A_0\|_{H^{2,p}} \leq  c(\Sigma, p)\|\langle i\Phi , D_0 \Phi\rangle\|_{L^p} .
$$
For $p<2$,
\begin{align}
\|\langle i\Phi , D_0\Phi\rangle\|_{L^p}  
  &\leq \|\Phi\|_{L^{\frac{2p}{2-p}}} \|D_0\Phi\|_{L^2}  \\
  & \leq c(\mcn_0 , \Sigma)\|D_0 \Phi\|_{L^2}.
\end{align}
Therefore, together with \eqref{estao}
\be
\|A_0\|^{2}_{H^{2,p}} +\|dA_0\|^{2}_{L^2} \ \leq \ c \|D_0 \Phi\|^{2}_{L^2}.
\ee
Therefore,
\be\la{estaoint}
\int_{0}^{\infty} \big( \|A_0\|^{2}_{H^{2,p}} +
\|dA_0\|^{2}_{L^2}\big) e^{2\delta t} \ dt \ \leq \ \infty   
\ee
by \eqref{dophi}.
By the equation for $\Phi$ in  \eqref{bom} we also have
\begin{align}
\|\pt\Phi\|_{L^2} \ 
  &\leq \ \|A_0\Phi\|_{L^2} + c\left(\|\eta\|_{\honea}
+\|v\|_{H^1} \right)\\
  & \leq c(\mcn_0,\Sigma) \left(\|D_0 \Phi\|_{L^2} + \|\eta\|_{\honea} +
\|v\|_{\hone}\right) \ 
\end{align}
and so
\be\la{estphiint}
\int_{0}^{\infty} \|\pt\Phi\|^{2}_{L^2} \, e^{2\delta t} \ dt \ \leq
\ \infty  \  .
\ee
Similarly for $\A$, we obtain
\be\label{estaint}
\int_{0}^{\infty} \|\pt \A\|^{2}_{L^2}\, e^{2\delta t} \ dt \ \leq
\ \infty  \  .
\ee
These estimates immediately imply \eqref{dot} by the equations
\eqref{bom} and then by \eqref{int} we conclude  that 
$$
\pt ({\mathbf A},\Phi) =  f  \mbox{ on } [0,\infty)\times \Sigma
$$
where $f\in L^2$ and $\int_{0}^{\infty} \|f\|_{L^2}(t) e^{2\delta t}
dt <\infty$, and thus
$$
({\mathbf A}_\infty, \Phi_\infty)\equiv lim_{t\to\infty}( {\mathbf A}(t),\Phi(t))
$$ 
exists in $L^2$ with  exponential rate of convergence. 

From these estimates and the above we have that $\Phi\in H^1$ and so is in $L^p$ 
for all $p<\infty$ and as $t\to\infty$, $\Phi\to\Phi_\infty$ strongly in $L^2$, 
$D_\A\Phi\to D_{\A_\infty}\Phi_{\infty}$ weakly in $L^2$, hence 
$\Phi D_\A\Phi  \to \Phi_\infty D_{\A_\infty}\Phi_\infty $ 
weakly in $L^1$;  in addition, as $\Phi$ is bounded in $L^q$ for any
$q<\infty$  (from $\Phi \in H^1$), the product $\Phi D_\A\Phi$  is bounded in $L^q$ for every $q<2$. This implies weak convergence 
in $L^q$ for $q<2$ so that $div (i\Phi, D_\A\Phi)$ converges weakly in $H^{-1, q}$; by 
Calderon-Zygmund then $A_0$ converges to a limit $A_{0\infty}$ weakly
in $H^{1,q}$ for each $q<2$ and so strongly in every $L^p$ for
$p<\infty$.   
Clearly by \eqref{estaoint} this limit is zero.

This completes the proof of the theorem in the case $|\Sigma| > 4\pi N$. 
\qed
  
\subsection{ Proof in the case $|\Sigma|<4\pi N$}

In this case a different approach is needed because 
the condition $\|\Phi\|^{2}_{L^2}\geq m>0$ which held and
was  used in the previous section is no longer valid (as is obvious
from the explicit form of the lower bound in \eqref{vpsn}).  Indeed,
as first observed in \cite{brad88}, in this case there are
no solutions to the Bogomolny equations \eqref{bogo}. To see this fact,
integrate the equation $v=0$ over $\Sigma$ to deduce
\beq
2\pi N-\frac{|\Sigma|}{2}+\frac{\|\Phi\|_{L^2}^2}{2}=0
\eeq
which is an impossibility when $4\pi N|>|\Sigma|$. We will
show instead that the minimizers of $\mcn$ 
have $v=\frac{l}{|\Sigma|}=constant$ and $\|\Phi\|_{L^2}= 0$, i.e. the field 
$\Phi$ is identically zero and $D_{\mathbf A}=\nabla-i\A$ is a constant 
curvature connection on $L$. Notice that the value of 
the constant $l$ is fixed as
\beq
l=2\pi N-\frac{|\Sigma|}{2}>0,
\eeq
since $\int_\Sigma vd\mu_g=l+\frac{1}{2}\|\Phi\|_{L^2}^2$ by
integration
of \eqref{bogo2}. Now
define $d>0$ by $d^2=|\Sigma|/\int_{\Sigma}v^2 d\mu_g$, then
\begin{align}
2\Bigl(\int_{\Sigma}v^2d\mu_g\Bigr)^{\frac{1}{2}}
|\Sigma|^{\frac{1}{2}}=
d\int_\Sigma v^2d\mu_g+\frac{|\Sigma|}{d}&=
2\int_{\Sigma}vd\mu_g+\int_{\Sigma}\bigl(
\sqrt{d}v-{\frac{1}{\sqrt{d}}}
\bigr)^2d\mu_g\notag\\
&=2(l+\frac{1}{2}\|\Phi\|_{L^2}^2)+\|(\sqrt{d}v-{\frac{1}{\sqrt{d}}})\|^2_{L^2}.
\notag
\end{align}
Referring to \eqref{bogdec},
this implies that the energy $\mcn$ can be written
\beq\la{bogdec2}
\mcn({\mathbf A},\Phi)
= 2\|\barpa \Phi\|_{L^2}^2+
\frac{1}{8|\Sigma|}\Bigl[
2(l+\frac{1}{2}\|\Phi\|_{L^2}^2)+\|(\sqrt{d}v-{\frac{1}{\sqrt{d}}})\|^2_{L^2}
\Bigr]^2\,+\,\pi N.
\eeq

From this it follows, since $l>0$, that
$\mcn\geq \pi N+\frac{l^2}{2|\Sigma|}$ and that this
lower bound is achieved with $\Phi$ identically zero
and $B=\frac{2\pi N}{|\Sigma|}=constant$, so that
$v=\frac{l}{|\Sigma|}=\frac{1}{d}=constant$. Thus
although the Bogomolny bound $\pi N$ is not
itself achieved in the case 
$4\pi N>{|\Sigma|}$ we have identified the greatest lower
bound and shown that it is achieved
with the ``pure magnetic'' constant curvature connections:
\begin{lemma}\la{min}
For $|\Sigma|<4\pi N$ 
$$
\min_{H^1(L\oplus\Omega^1)}\mcn=
\pi N+\frac{(4\pi N-{|\Sigma|})^2}{8|\Sigma|}$$
and 
this minimum value is achieved by taking the Higgs field $\Phi$
identically zero and $D$ to be a constant curvature connection.
\end{lemma}

To analyze the gradient flow in this case it is useful to introduce the
variable $y=v-l$ (in place of $v$), so that \eqref{ebm}
are replaced by
\begin{align}\la{ebm2}
\begin{split}
&D_0\eta-4\barpa(e^{-2\rho}\pbfa\eta)+(l+|\Phi|^2)\eta
=-y\eta,\\
&(\pt-\Delta+|\Phi|^2)y=-|\Phi|^2l-4|\eta|^2,\\
&D_0\Phi-4e^{-2\rho}\pbfa\barpbfa\Phi+l\Phi=-y\Phi.\\
\end{split}
\end{align}

Corresponding to \eqref{enid} we have the following
integral inequalities for solutions of \eqref{ebm2}

\begin{align}
\begin{split}\la{enid2}
&e^{2\delta T}\|\eta(T)\|_{L^2}^2+2\int_0^Te^{2\delta t}\Bigl(
Q_{\A,\Phi}(\eta)+(l-\delta)\|\eta\|_{L^2}^2
\Bigr)\,dt\;\leq\;\|\eta(0)\|_{L^2}^2-2\int_0^T\int_\Sigma e^{2\delta t}y|\eta|^2 
e^{-2\rho}d\mu_g dt,\\
&e^{2\delta T}\|y(T)\|_{L^2}^2+2\int_0^Te^{2\delta t}\Bigl(
Q_{\Phi}(y)-\delta\|y\|_{L^2}^2
\Bigr)\,dt\;\leq\;\|y(0)\|_{L^2}^2
-2\int_0^T\int_\Sigma e^{2\delta t}y(|\Phi|^2l+4|\eta|^2)
e^{-2\rho}d\mu_g dt,\\
&e^{2\delta T}\|\Phi(T)\|_{L^2}^2+2\int_0^Te^{2\delta t}\Bigl(
\|\barpbfa\Phi\|_{L^2}^2+(l-\delta)\|\Phi\|_{L^2}^2
\Bigr)\,dt\;\leq\;\|\Phi(0)\|_{L^2}^2
-2\int_0^T\int_\Sigma e^{2\delta t}\,y|\Phi|^2\,
e^{-2\rho}\,d\mu_g dt.
\end{split}
\end{align}

Notice that since $l>0$, there is a natural mechanism
forcing $\Phi$ to converge to zero at an exponential rate.
However this fact necessitates modification of 
the arguments based
on the lower bounds in lemma \ref{lb}, which are dependent
upon $\|\Phi\|_{L^2}\geq m>0$. The presence of $l>0$ in the
integrals in the first and third identities  means that for
$\delta<\frac{l}{2}$ there exists $\tilde\gamma=\tilde\gamma
(\mcn_0,l)$ such that
\begin{align}
Q_{\A,\Phi}(\eta)
+(l-\delta)\|\eta\|_{L^2}^2&\geq\tilde\gamma\|\eta\|_{\honea}^2,
\quad\hbox{and}
\\
\|\barpbfa\Phi\|_{L^2}^2+(l-\delta)\|\Phi\|_{L^2}^2
&\geq\tilde\gamma\|\Phi\|_{\honea}^2.
\end{align}
For the middle identity, recall that if $\overline y=
\frac{1}{|\Sigma|}\int_\Sigma y$ then Poincare's inequality
says that 
$$
\int e^{-2\rho}\sum_{k}|\partial_ky|^2d\mu_g \geq C_P\|y-\overline y\|_{L^2}^2
\geq C_P(\|y\|_{L^2}^2-\|\overline y\|_{L^2}^2),$$ so that
there exists $\gamma'>0$ such that
$$
Q_\Phi(y)\geq\gamma'(\|y\|_{H^1}^2-
\|\overline y\|_{L^2}^2).
$$
But also $l+\overline y=\frac{1}{|\Sigma|}\int_\Sigma v=
l+\frac{1}{2|\Sigma|}\|\Phi\|_{L^2}^2$, so that
$\overline y=\frac{1}{2|\Sigma|}\|\Phi\|_{L^2}^2$. Therefore,
overall we have the following set of controlling inequalities:

\begin{align}
\begin{split}\notag
&e^{2\delta T}\|\eta(T)\|_{L^2}^2+2\int_0^Te^{2\delta t}
\tilde\gamma\|\eta\|_{\honea}^2
\,dt\;\leq\;\|\eta(0)\|_{L^2}^2-2\int_0^T\int_\Sigma 
e^{2\delta t}y|\eta|^2 
e^{-2\rho}d\mu_g dt,\\
&e^{2\delta T}\|y(T)\|_{L^2}^2+2\int_0^Te^{2\delta t}
\Bigl((\gamma'-\delta)\|y\|_{H^1}^2-\gamma'\|\overline y\|_{L^2}^2
\Bigr)\,dt\;\;\;\;\leq\;\|y(0)\|_{L^2}^2\\
&\phantom{e^{2\delta T}\|y(T)\|_{L^2}^2+2\int_0^Te^{2\delta t}
\Bigl((\gamma'-\delta)\|y\|_{H^1}^2-\gamma'\|\overline y\|_{L^2}^2
\Bigr)\,dt\;\;\;\;\leq\;}
-2\int_0^T\int_\Sigma e^{2\delta t}y(|\Phi|^2l+4|\eta|^2)
e^{-2\rho}d\mu_g dt,\\
&e^{2\delta T}\|\Phi(T)\|_{L^2}^2+2\int_0^Te^{2\delta t}
\tilde\gamma\|\Phi\|_{\honea}^2\,dt\;\;\;\;
\leq\;\|\Phi(0)\|_{L^2}^2
-2\int_0^T\int_\Sigma e^{2\delta t}y|\Phi|^2
\,e^{-2\rho}\,d\mu_g\, dt
\end{split}
\end{align}
Add the three inequalities, and use the fact that
$\|\overline y\|_{L^2}^2=|\Sigma||\overline y|^2=
\frac{1}{2}\overline y\|\Phi\|_{L^2}^2$ to absorb the
negative term in the integral in the second inequality
by the integral in the third one. Next bound the
nonlinear terms on the right hand sides in the same 
way as in \eqref{lineqb1}-\eqref{int}, and conclude that as long as
$\|y(0)\|_{L^2}$ is initially small then
\begin{align}
\begin{split}\notag
&e^{2\delta T}\|\eta(T)\|_{L^2}^2
\leq\;\|\eta(0)\|_{L^2}^2,\\
&e^{2\delta T}\|y(T)\|_{L^2}^2\leq\;\|y(0)\|_{L^2}^2,\\
&e^{2\delta T}\|\Phi(T)\|_{L^2}^2
\leq\;\|\Phi(0)\|_{L^2}^2,
\end{split}
\end{align}
holds for all $T>0$, and
\beq\la{mmm}
\int_0^\infty e^{2\delta t}\Bigl[
\|y\|_{H^1}^2
+\|\Phi\|_{\honea}^2
+\|\eta\|_{\honea}^2\Bigr]
\,dt\, <\infty
\eeq
from which convergence can be deduced as in the previous
section: \eqref{dophi} and \eqref{estaoint} hold as a consequence of
\eqref{mmm}, and hence by the first two equations
of \eqref{bom} the estimate \eqref{estaint} also holds, and so 
$\A(t)$ converges to a limit $\A_\infty$ at an exponential rate.
The proof of the main theorem is now complete.\qed

\end{document}